\documentclass [12pt] {report}
\usepackage[dvips]{graphicx}
\usepackage{epsfig, latexsym,graphicx}
\usepackage{amssymb}
\newcommand {\D}[2] {\displaystyle\frac{\partial{#1}}{\partial{#2}}}

\newcommand {\Dd}[3] {\displaystyle\frac{\partial^2{#1}}{\partial{#2}\partial{#3}}}

\newcommand {\ga} {\gamma}

\newcommand {\de} {\delta}
\newcommand {\prtl} {\partial}
\newcommand {\fr} {\displaystyle\frac}

\newcommand {\be} {\begin{equation}}
\newcommand {\ee} {\end{equation}}
\newcommand {\ba} {\begin{array}}
\newcommand {\ea} {\end{array}}
\newcommand {\bp} {\begin{picture}}
\newcommand {\ep} {\end{picture}}
\newcommand {\bc} {\begin{center}}
\newcommand {\ec} {\end{center}}
\newcommand {\bt} {\begin{tabular}}
\newcommand {\et} {\end{tabular}}
\newcommand {\lf} {\left}
\newcommand {\rg} {\right}

\newcommand {\cC} {{\cal C}}
\newcommand {\cF} {{\cal F}}
\newcommand {\cH} {{\cal H}}
\newcommand {\cI} {{\cal I}}
\newcommand {\cM} {{\cal M}}

\newcommand {\cR} {{\cal R}}
\newcommand {\cS} {{\cal S}}

\newcommand {\ses} {\medskip}

\newcommand {\bibit} {\bibitem}
\newcommand {\nin} {\noindent}

\newcommand {\cU} {{\cal U}}   \newcommand {\cY} {{\cal Y}}
\newcommand {\cV} {{\cal V}}

\newcommand {\Cos} {\mathop{\rm Cos}\nolimits}

\newcommand {\Sin} {\mathop{\rm Sin}\nolimits}
\newcommand {\Tan} {\mathop{\rm Tan}\nolimits}

 \usepackage{amsmath}

\def\2#1#2#3{{#1}_{#2}\hspace{0pt}^{#3}}
\def\3#1#2#3#4{{#1}_{#2}\hspace{0pt}^{#3}\hspace{0pt}_{#4}}
\newcounter{sctn}
\def\sec#1.#2\par{\setcounter{sctn}{#1}\setcounter{equation}{0}
                  \noindent{\bf\boldmath#1.#2}\bigskip\par}

    \newcommand {\s} {\ses}

\newcommand{\ed}{\end{document}}

    \newcommand{\hU}{\hat U}   \newcommand{\hZ}{\hat Z}

\newcommand{\ccUss}{U^{\{c\}} }

\newcommand{\ccfss}{f^{\{c\}} }   

\newcommand{\bvU}{\breve U }   
   \newcommand{\bvZ}{\breve Z }   
\newcommand{\bvf}{\breve f }  
    
  \newcommand{\bvt}{\breve t }

   \newcommand{\bvth}{\breve{\theta} }

\newcommand{\tht}{\theta }   

    \newcommand{\rmsps}{Riemannian  space}

  \newcommand{\rmt}{Riemannian  metric tensor }    

\newcommand{\rmct}{Riemannian  curvature tensor }

 \newcommand{\frn}{Finslerian }

 \newcommand{\ffss}{Finsleroid space}

\newcommand{\fsp}{Finsler space } 
 \newcommand{\ffsps}{Finsleroid space}

  \newcommand{\FMF}{Finslerian metric function }   \newcommand{\FMFs}{Finslerian metric function}
  \newcommand{\fmf}{Finslerian metric function }   \newcommand{\fmfs}{Finslerian metric function}
  \newcommand{\fmt}{Finslerian metric tensor }    

  \newcommand{\imt}{indicatrix metric tensor }      

\newcommand{\frd}{Finsleroid }   \newcommand{\frds}{Finsleroid}

\newcommand{\FRD}{Finsleroid }

\newcommand{\frdmt} {Finsleroid metric tensor }

\newcommand{\MT} {metric tensor }

\newcommand{\FMT}{Finslerian  metric tensor } \newcommand{\FMTt}{Finslerian  metric tensor}

\newcommand{\ffmf}{Finsleroid metric function }   \newcommand{\ffmfs}{Finsleroid metric function}

 \newcommand{\ffmts}{Finsleroid metric tensor}
\newcommand{\ffind}{Finsleroid indicatrix }

 \newcommand{\IND}{indicatrix }    

 \newcommand{\ind}{indicatrix }    \newcommand{\inds}{indicatrix}

  \newcommand{\indsp}{indicatrix space }      

\newcommand{\indmt}{indicatrix  metric tensor }

  \newcommand{\irsp}{indicatrix Riemannian space }  

    \newcommand{\rsp}{Riemannian space }

\newcommand{\IRSR}  { $\cR^{ Ind}_x$ }  

\newcommand{\bIRSR}  { $\breve\cR^{ Ind}_x$ }

\newcommand{\htph}{\{H;\theta,\phi\}}

    \newcommand{\rhpt}{ right-hand part}

\newcommand{\CT}{ Conclusive Theorem }

\newcommand {\SSch} {Separation Scheme }   \newcommand {\SScht} {Separation Scheme}

\newcommand{\FH}{${\cal FH}_x^{P_3}$ }

 \newcommand{\UZ}{${\cal UZ}$}
\newcommand{\UZfsp} {${\cal UZ}$-Finsler space }  \newcommand{\UZfsps} {${\cal UZ}$-Finsler space}

   \newcommand{\FUZs}{$F^{\cal UZ}$}
\newcommand{\FE}{$\cal FRD^{{UZ}} $}

\newcommand {\hUZ}{$h^{\cal UZ}$ } 
\newcommand {\hUZss}{h^{\cal UZ} }

\newcommand {\gUZ}{$g^{\cal UZ}$ } 
\newcommand {\gUZss}{g^{\cal UZ} }

\newcommand {\DUZss}{D^{\cal UZ} }

\newcommand{\calcns}{calculations }

\newcommand{\Prpn}{Proposition }  

\newcommand{\reprn}{representation }  
\newcommand{\reprns}{representations }

  \newcommand{\defns}{definition}

\newcommand{\hty}{homogeneity }  

\begin {document}


{

{



\s\s

\s\s


\s


\s



\newcommand {\vsec}{vertical  section }   

\s


\s

\newcommand{\vsecFRD} {$\cV$-section of the \frd }   

\s \s  \s

        
        

\s

        


\s

        
        

\s

        
        
        

\s


        
        

\s

        
        

\s







\s\s


\newcommand{\HxPss}{ {\cal H}_{x;P_3} }

\s

\newcommand{\RHx}{$R^{P_3}_x$ } \newcommand{\RHxss}{R^{P_3}_x }

 \newcommand{\zHxss}{z^{P_3}_x }

\newcommand{\thtHx}{$\theta^{P_3}_x$ } \newcommand{\thtHxss}{\theta^{P_3}_x }


\s\s


        




\s

\newcommand {\hsec}{horizontal  section }   

\s

   \newcommand{\FHxPss}  { {\cal FH}_{x;P_3}}

\s

\newcommand{\hsecFRD} {$\cH$-section of the \frd }   

\s\s

\s \s  \s

        
        

\s

        
        

\s


        
        

\s

        
        
\s

\s


\s  \s\s

        
        

\s \s\s


\s\s\s

        
        










\newcommand{\azmanglth}{azimuthal angle $\tht$ }


  \newcommand{ \chfns} {characteristic function}

\newcommand{\strln} {straight line }   \newcommand{\strlns} {straight line}





\begin{titlepage}

\vspace{0.1in}

\begin{center}

{\bf \Large Finsleroids with three axes in  dimension $N=3$
}

\end{center}

\vspace{0.3in}

\bc

\vspace{.15in}{\large G.S. Asanov\\}  \vspace{.25in}
{\it Devision of Theoretical Physics, Moscow State University\\
119992 Moscow, Russia\\
{\rm (}e-mail: asanov@newmail.ru{\rm)}} \vspace{.15in}

\ec

\begin{abstract}

\s  \s

The \CT has been  established to determine  the dependence of
the   three-axes positive-definite {\ffmfs}s
 $F$
on the \frd \azmanglth
in the three-dimensional case $N=3$, provided that
 the condition of  the angle-separation in the  involved characteristic functions is implied.
The complete   set
  of algebraic and  differential equations
is derived in all rigor which
 are necessary and sufficient in order that the  function  $F$ belong to the class.
 It proves possible to solve the equations and obtain the explicit dependence of the involved
 characteristic functions on the angle
 $\tht$.

\s

\nin

\s

\nin

\end{abstract}

\end{titlepage}

\vskip 1cm



   
\setcounter{equation}{0}

\s \s



\nin
{\bf  \large 1   Introduction and preliminaries}



\s \s

\nin
In the preceding publications [1-3] the method  was proposed and used to derive the class of
the two-axis pseudo-Finsleroids under the structural pre-assumption of the angle separation
in the characteristic functions. The consideration was four-dimensional, $N=4$, and the signature of \FMT was
                  presumptively  (+---).
In the present paper the similar method is  applied to find the three-axis Finsleroids in the positive-definite
three-dimensional case, so that  $N=3$ and the signature is $(+++)$.
The treatment as well as the evaluation process are systematic.


The consideration uses  the property of separation of dependence  with respect to the \frd angles $\theta$ and $\phi$.
The evaluation option of the property is highly significant, namely  the  separation  makes it possible to reduce
 the  characteristic  partial  differential equations to a set of the ordinary differential equations with respect to $\theta$ and $\phi$.
As we shall show below,  the equations can be solved to give us the explicit algebraic $\tht$-angle representations
for all the involved characteristic   \frd
functions.

The origin of the  Finsler Geometry  can be traced back to the famous publications [4-6].
During the succeeding decades the geometry has been intensively elaborated  in many  interesting and important directions (see [7-13]). At present, the Euclidean metric  yields  the profound geometric base for various scientific  theories.
The metric is constructed with the help of the quadratic form of tangent vectors to meet the simplest algebraic patterns. To go over the simplest    presumption we can naturally apply the methods of the Finsler geometry.

Obviously, the applied ability of the Finslerian methods
 is proportional to the variety of the particular Finslerian
metric functions which can be proposed and well-developed to geometrically underline the possible anisotropic distributions of vectors in the tangent spaces of differentiable manifolds.

It is interesting and constructive to retain   the condition of constancy of the curvature of indicatrix, in which case we qualify the Finsler metric function to be of the {\it Finsleroid} type in the positive-definite case.
The one-axis Finsleroid metric function  reveals various interesting properties (see [14,15] and references therein) of the fundamental significance to develop applications.
In the present paper, we confine attention to the Finsleroid metrics  possessing three axes of anisotropy.
Our analysis is systematic, the presentation is concise and complete.


A direction in the tangent space is said to be {\it geometrically distinguished} if the direction makes trace in the structure
of the \fmfs.
 We are aiming in the present paper to solve the following problem in the three-dimensional case.

\s\s

\nin
{\bf Problem} Find the
{\fmfs}s
 $F$ which entail the angle representation (1.3)   of the
associated \fmt  $g$
in the case when the functions $F$ are assigned by three geometrically distinguished axes, provided that the property of separation of the angle dependence is fulfilled.

\s\s



\nin
To this end, we consider  a three-dimensional differentiable manifold $M$ and treat the tangent bundle $TM$ as the unification of the tangent spaces
 $T_xM\subset TM, x\in M$. The tangent vectors  $y\in T_xM$ are supported by the points $x\in M$.
Let $F$ be a  \FMFs, such that $F=F(x,y)$ and the positively first-degree homogeneity
with respect to the vector argument $y$ together with sufficient smoothness is implied.
The {\it\ind } ${\cI}_xM\subset T_xM$ in a tangent space $T_xM$ is the surface   defined by
  $${\cal I}_xM=\{y\in T_xM: F(x,y)=1\}.$$
 The definition extends the notion of the Euclidean sphere.



In each tangent space $T_xM$ the \FMT induces a \rmt
on the \ind
(according to the definition (2.1)),
so that the indicatrix  ${\cI}_xM$  becomes a \rmsps,
 to be denoted by
 $\cR_x^{Ind}$.

 \s \s

\nin
{\bf Definition 1.1} If at a point $x\in M$ of support   the \IND
 space  $\cR_x^{Ind}\subset T_xM $ is of constant positive curvature, then the \IND body
  \be
  \cF_x=\{y\in T_x: F(x,y)\le1\}
  \ee
is called the {\it \frds}.

\s\s

\nin  This definition extends the notion of the Euclidean ball.


  If the  space  $\cR_x^{Ind}\subset T_xM $ is of constant positive curvature at a point
$x\in M$ of support, we call $F$ the {\it \ffmf}  at this point.
If this property fulfills at any point $x\in M$ of the space,  we call $F$ the {\it \ffmf on $TM$}.
We shall denote the function $F$ by  $F^{Frd}$
to emphasize the \frd nature of the function.

The  \FMT $g$ constructed from the  \ffmf $F^{Frd}$
 will be called the {\it \ffmts}.


Accordingly, we shall use the following definition.

\s \s

\nin
{\bf Definition 1.2} The space
\be
\cF_{\{3\}}=\{M;TM;y\in TM;F^{Frd};g\}
\ee
is called  the three-dimensional
 {\it \ffss.}


\s  \s


\nin
Since   we confine the consideration to the  three-dimensional case,  the tensorial indices $i,j,...$
will take on the values (1,2,3).
 Our consideration will be of local nature,   so that we shall often
represent the vectors and tensors by means of their local components with respect to the natural frame,
in particular
 $y=\{y^i\}$.

Following the known methods of the Finsler geometry, we  construct the tensor
$g=g(x,y)$ by means of the derivative components:
$ g_{ij}=(1/2)\partial^2 F^2/\prtl y^i/\prtl y^j$ and
 $g=\{g_{ij}\}$.
We shall use also the tensor $\{h_{ij}\}$ obtained from the expansion $g_{ij}=l_il_j+h_{ij}$,
where $l_i=\partial F/\prtl y^i$ are covariant components of the  unit vector $l=\{l^i\}$, such that
$l^i=y^i/F$, $l_i=g_{ij}l^j$, and $F(x,l)=1.$


\s

  Let us denote by  $\theta=\theta(x,y)$ and $\phi=\phi(x,y)$ the  scalars  on the manifold $M$
  assuming that  the functions $\theta(x,y)$  and $\phi(x,y)$ are functionally independent  and smooth of at least class $C^2$
with respect to the variable $y$ in each tangent space $T_x$,
and also  are positively homogeneous of the  degree zero with respect to  vectors $y\in T_x$.
We assume that the Finsler space possesses the following property:
two   scalars $\theta$ and $\phi$ exist such that the \FMT of the space can be written in the form
\be
g_{ij}=l_il_j+\fr1{\cal C}\lf(\theta_i\theta_j+\sin^2\theta \phi_i\phi_j\rg)F^2, ~~  {\cal C}={\cal C}(x)>0,
\ee
where
the notation $\{\theta_i=\prtl\theta/\prtl y^i, ~ \phi_i=\prtl\phi/\prtl y^i\}$
has been used.


 \s


In Sect. 2, the important aspects of the \ind geometry required for our \frd study are described.

In Sect. 3, we are inspiring by the attractive idea
to construct the   \FRD space with three axes  on a three-dimensional  \rsp
 endowed with three vector fields. Such a \rsp
  is
 the geometrical background for the theory developed in the present paper.
 The involved vectors are assumed naturally to be  linearly independent
 at each point of the background manifold. They are also used to define three axes for the \frds s.
 In terms of these axes, we separate the angle dependence
of the  characteristic functions
which enter the  functions $F$ in accordance with the
 Separation Scheme (3.9).
The respective notion of the
\UZfsp
 is introduced  such  that the \fmf $F$ of the space  gets compatible with
the Separation Scheme.
The interesting concepts  of the  \hsec and the \vsec
of the \frd are appeared  to study.

We also
 observe the attractive possibility to extend the Euclidean trigonometric functions
to obtain the generalized  \FRD trigonometric functions
retaining the Euclidean structure of dependence on the angles $\theta$ and $\phi$.
All the components of the
\UZfsp
 tensor $\hUZss_{ij}=\gUZss_{ij}-l_il_j$ are evaluated in the concise and explicit form.
The entailed determinant
$\DUZss=\det(\gUZss_{ij})$
proves to be  the product $\DUZss=D_1D_2$
possessing  the property of separation of dependence on  $\tht$ and $\phi$, namely
$D_1=D_1(x,\theta)$ and $ D_2=D_2(x,\phi).$

\s


In Sect. 4,
we clarify when
 the  \UZ-Finsler space is of the \FRD type.
To this end
 we compare  the obtained components of the tensor
   $\hUZss_{ij}$
 with the respective components of the \frds-tensor $h_{ij}^{\{H;\theta,\phi\}}$
written in (2.7). It is surprising enough, that this method which is simple conceptually (although includes rather lengthy
calculations) leads to Conclusive Theorem which formulates the general result with the help of two rather simple differential equations
 (4.3)-(4.4) which we call the
{\it \UZ-Finsleroid characteristic equations}, the
{\FE-equations} for short.
It proves possible to solve the equations to clarify the
dependence of the characteristic functions $U$ and $f$ on the angle $\theta$. The solutions are explicitly found in terms of simple   algebraic functions (see (4.19) and (4.22)).
The  \UZ-Finsler space is the \FRD space at any choice of the dependence of the involved function $Z$ on $t$, provided
that the generating function $U$ fulfills the \FE-equations.

In the last section we emphasize
the most important aspects of the developed theory.




\setcounter{equation}{0}

\s\s

\nin
{\bf \large 2 Differential geometry of \frd \ind}

\s\s

\nin
The following definition gives rise to  significant geometrical ideas and consequences
in our analysis.

\s\s

\nin
{\bf
Definition 2.1} Formula  (1.3) introduces the {\it angle representation}  of the  \FMTt,
interpreting $\tht$ and $\phi$ as the angle variables.

\s\s

\nin
Indeed,
let us denote $u^1=\theta$ and $u^2=\phi$, fix a point $x\in M$, and interpret the set
 $\{u^a\}$ to be a coordinate system on the indicatrix supported by $x\in M$. The indices $ a,b,...$ will be specified on the range
  (1,2).
  Parameterize the unit vectors $l=\{l^i\}$  of the tangent space $T_x$ by the help of the
coordinates $u^a$
  using the respective vector field $t^i=t^i(x,u)$. We get the parametrical representation
$l^i=t^i(x,u)$ of unit vectors, which in turn gives rise to  the so-called {\it projection factors} $t^i_a=\partial t^i/\partial u^a$ (see Sect. V.8 in [7]).
With these objects, we  can construct the \MT   $\breve i=\{i_{ab}\}$ on the \ind  by means of the components $i_{ab}=g_{ij}t^i_at^j_b$. Owing to  the identities $l_it^i_a=0$ and the expansion $g_{ij}=l_il_j+h_{ij}$, we can write simply
\be
 i_{ab}=h_{ij}t^i_at^j_b.
\ee

\s

This  equality  can be inverted with the help of the vectors $u^a_i=\partial u^a/\partial y^i$. Indeed,
upon contracting the equality  with $u^a_iu^b_j$ and taking into account the  \frn  identity $Ft^i_au^a_j=h^i_j$, we get
the  inverse  \reprn
\be
h_{ij}=i_{ab}u^a_iu^b_jF^2.
\ee
This is a general \frn result.


The  tensor $\breve i=\breve i(x,u)$ thus obtained is called the {\it  induced \rmt  on \inds},
the {\it \indmt} for short.
 From Formulas (1.3) and (2.1) we get
\be
i_{11}=\fr1{\cC}, ~ i_{12}=0, ~ i_{22}=\fr1{\cC }\sin^2\theta.
\ee
Using these components, we can  construct the
\rmct
 $R^f{}_{abd}=
R^f{}_{abd}(x,u)$
and  directly obtain the \reprn
\be
 R_{abcd}=  \cC  (i_{ad}i_{bc}-i_{ac}i_{bd}), ~~     R_{eabd} = i_{ef}R^f{}_{abd},
\ee
 which just  tells us that the \ind   is a space of the constant curvature $\cC$.


 The verification of the \reprn (2.4) is a short task. Indeed, constructing the Christoffel symbols
$
i^c{}_{ab}=i^{ce}i_{e,ab},   ~~  i_{e,ab}=(i_{ae,b}+i_{be,a}-i_{ab,e})/2,
$
we obtain the components
$
i^1{}_{11}=  i^1{}_{12}= i^2{}_{11}= i^2{}_{22}=0, ~~
i^1{}_{22}=-\sin\theta    \cos\theta, ~~
i^2{}_{12}=     \cos\theta/\sin\theta.
$
 Evaluating the indicatrix curvature tensor
\be
R^f{}_{abd} = \D{i^f{}_{ab}}{u^d}-\D{\ga^f{}_{ad}}{u^b}
+i^e{}_{ab} i^f{}_{ed}
-i^e{}_{ad} i^f{}_{eb}
\ee
and then lowering the index to obtain the tensor
$ R_{eabd} = i_{ef}R^f{}_{abd} $,
we can observe that
$$
\cC R_{1212} =   \D{i^1{}_{21}}{u^2}-\D{i^1{}_{22}}{u^1}
+i^e{}_{21} i^1{}_{e2}
-i^e{}_{22} i^1{}_{e1},
  $$
or
$$
\cC R_{1212} =   -\D{i^1{}_{22}}{u^1}
+i^2{}_{21} i^1{}_{22}=
-\sin^2\theta +   \cos^2\theta  -\cos^2\theta
=-\sin^2\theta.
  $$
So the equality $ R_{abcd}=  \cC  (i_{ad}i_{bc}-i_{ac}i_{bd})$  is true.

\s\s

Thus the following assertion  is valid.

\s \s

\nin
{\bf Primary  Theorem 2.1} {\it If the \fmt  of a  \fsp
 admits the  \reprn \rm(1.3),
\it then the \fsp
 is of the \frd
  type,
namely the  \indsp  \IRSR is a space of the constant curvature   ${\cal C}(x)$.}


\s\s

\nin
Geometrically, our constructions in the present paper are founded by the following \defns.
\s\s

\nin
{\bf Definition 2.2}  The space    \bIRSR=$\{\cI_xM, \breve i_x, u^a\}$
is called the {\it \frd \irsp} supported by the point $x\in M$, where the \imt
$\breve i_x$ has the components $\{i_{11},i_{12},i_{22}\}$ listed in (2.3).

\s\s

\nin
We introduce the scalar $H=H(x)$ which will play the role of the {\it  characteristic Finsleroid parameter}
in the \ffsps.
The geometrical  meaning of the scalar is that the square  $H^2=(H(x))^2$  taken at a point $x$ is equal to the value of curvature of the
  \ffind
constructed in $T_x$.
Accordingly, we make the substitution ${\cal C}=H^2 $.
The value of $H$ will be restricted by the condition
$ 0<H\le 1 $ implied to hold at any admissible point $x\in M$.
The premised \reprn  (1.3) of the metric tensor can conveniently be written in the form

\be
g_{ij}=l_il_j+h_{ij}^{\{H;\theta,\phi\}}
\ee
with
\be
h_{ij}^{\{H;\theta,\phi\}}=\fr1{H^2}\lf(\theta_i\theta_j+\sin^2\theta \phi_i\phi_j\rg)F^2.
\ee


\s

Because of  (2.3),
the squared length element $(ds)^2=i_{ab}du^adu^b$ on the \frd\ind reads
 \be
(ds^2)|_{Frd.}=  \fr1{H^2}\lf((d\theta)^2+\sin^2\theta(d\phi)^2\rg),
\ee
which extends the Euclidean \reprn
$$
(ds^2)|_{Eucl.}=  (d\theta)^2+\sin^2\theta(d\phi)^2
$$
by means of the factor $1/H^2$.


\s\s





\setcounter{equation}{0}

\nin
{\bf\large 3 Geometrical background and   separation of  angle dependence}

\s\s

\nin
The notion of the  \FRD space with three axes
 can naturally  be introduced on a three-dimensional  Riemannian space. Namely,
let $\cM^*_3$ be a three-dimensional differentiable manifold
which admits the introduction of three vector fields assuming that the fields are linearly independent at each point $x\in \cM^*_3$.
The associated tangent bundle $ T\cM^*_3$ is the unification of the respective tangent spaces $ T_x\cM^*_3$.
Denote by $i,j, i_{\{3\}}$ the one-forms which represent the vector fields.
 We shall apply
 the localized coordinate representations
$\{i=i_iy^i,j=j_iy^i, i_{\{3\}}=i_{\{3\}i}y^i\}$, where $i_i=i_i(x),j_i=j_i(x),$ and $ i_{\{3\}i}=i_{\{3\}i}(x)$,
such that $\{i_i,j_i,i_{\{3\}i}\}$ are the covariant vector fields on $\cM^*_3$ and $y=\{y^i\}$ denotes the tangent vector: $y\in T_x\cM^*_3$. Below, the short notation $T_x$
will be used instead of the complete $T_x\cM^*_3.  $


Constructing the  \rmt $a=\{a_{ij}(x)\}$ by means of the expansion
\be
a_{ij}=i_ii_j+j_ij_j+i_{\{3\}i}i_{\{3\}j},
\ee
we get the three-dimensional
Riemannian space  $\cR^*_3$ according to the following  definition.

\s  \s

\nin {\bf Definition 3.1} The notion
\be
\cR^*_3=\{\cM^*_3,i,j,i_{\{3\}},a\}
 \ee
 is called the {\it background  three-axes  \rmsps}.

\s \s

\nin
The next definition emphasizes the structural role of the introduced vectors $\{i_i,j_i,i_{\{3\}i}\}$.

 \s \s

\nin {\bf Definition 3.2}
 The members of the set $\{i_i,j_i,i_{\{3\}i}\}$ are called the {\it basis vectors} of the space $\cR^*_3$.

\s \s

\nin
The  indices of vectors and tensors in the space $\cR^*_3$ will be raised by means of the contravariant tensor  $\{a^{ij}(x)\}$
reciprocal to the covariant tensor $\{a_{ij}(x)\},$ such that
 $a_{in}a^{nj}=\de^j_i, j^i=a^{in}j_n, {\it etc;}$ $\de$ stands for the Kronecker symbol.



By  means of $\cS_x^{\{i\}}\subset T_x$, $\cS_x^{\{j\}}\subset T_x$, and  $\cS_x^{\{i_{\{3\}}\}}\subset T_x$
 we denote the          {\strlns}s
  which pass through the support center $O_x\in T_x$ of the tangent space $T_x$  and
 respectively include the basis vectors,  namely $i_x\in \cS_x^{\{i\}}$,  $j_x\in \cS_x^{\{j\}}$, and
 $i_{\{3\}x}\in \cS_x^{\{i_{\{3\}}\}}$. The inclusions
 $O_x\in \cS_x^{\{i\}}$,  $O_x\in \cS_x^{\{j\}}$, and
 $O_x\in \cS_x^{\{i_{\{3\}}\}}$ are implied.
Below we  consider the
{
 \frds s determined by the \FMFs s $F$ constructed by the help of these objects.    Accordingly, our treatment  complies with the following definitions.

 \s\s

\nin
{\bf Definition 3.3} The \strln
 $\cS_x^{\{i_{\{3\}}\}}$ is called
 the {\it vertical  Finsleroid  axis } at point $x\in \cM^*_3$.
The {\strlns}s
$\cS_x^{\{j\}}\subset T_x$ and  $\cS_x^{\{i\}}\subset T_x$ are called respectively
the {\it primary horizontal  Finsleroid  axis } and
the {\it secondary horizontal  Finsleroid  axis } at point $x\in\cM^*_3$.

\s \s

\nin
{\bf Definition 3.4} The space
\be
 \cF^*_3= \{\cM^*_3,i,j,i_{\{3\}},a,F,g\}
\ee
is called the
{\it   three-axes  Finsleroid space}. The entered $g$ denotes the
\fmt
 derived from the function $F$.

\s \s

\nin
In this way, we  construct the three-dimensional Finsleroid space $\cF^*_3$ on the three-dimensional
Riemannian space $\cR^*_3$.







\s\s




\nin
{\bf Definition 3.5} The plane
                                                   $\HxPss\subset T_x$
  which is orthogonal to the vertical axis
$\cS_x^{\{i_{\{3\}}\}}\subset T_x$  of the \frd in $T_x$ and intersects the axis at a point $P_3$ is called
the {\it \hsec of} the tangent space $T_x$.
The restriction  $\FHxPss  =\cF_x\cap \HxPss$
is called the
 {\it \hsec of the \frds,}
the {\it\hsecFRD} for short; $\cF_x$ is the \frd defined in (1.1).

\s\s




\nin
{\bf Definition 3.6} The plane $\cV_{x;\phi}\subset T_x$ which includes the vertical axis
$\cS_x^{\{i_{\{3\}}\}}\subset T_x$
of the \frd in $T_x$ and corresponds to a fixed  angle value $\phi$ is called
the {\it \vsec of}  $T_x$.
The restriction  ${\cal FV}_{x;\phi}=\cF_x\cap \cV_{x;\phi}$
is called the
 {\it \vsec of the \frds},
 the {\it {\it\vsecFRD}} for short.


\s  \s

\nin
{\bf Definition 3.7}
The space  $T^*_x\subset T_x$ which is obtained after deleting the vertical  Finsleroid  axis,
so that
$T^*_x=T_x  \setminus \cS_x^{\{i_{\{3\}}\}},$
is called the {\it  axially  reduced tangent space}. The respective bundle $T^*\cM^*_3\subset T\cM^*_3$ is called
the {\it axially  reduced  tangent bundle}.

\s  \s

\nin
{\bf Definition 3.8} The space
\be
\cF^{\{ \setminus \cS^{\{i_{\{3\}}\}}\} }_3=\{\cM^*_3,T^*\cM^*_3,y\in T^*_x\cM^*_3,i,j,i_{\{3\}},a,F,g\}
\ee
is called the
{\it axially reduced    three-axes   Finsleroid space.}

\s  \s


\nin
Using the basis vectors $\{i^i,j^i,i^{\{3\}i}\}$,
it is convenient to apply the expansion
\be
y^i=i^iy^1+j^iy^2+i_{\{3\}}^iy^3
\ee
in terms of the components
\be
\{y^a\}=\{y^1,y^2,y^3\}:  ~~  y^1=i, ~ y^2=j, ~ y^3=i_{\{3\}}.
\ee
Also, we introduce the notation
\be
 z=y^3\equiv i_{\{3\}}, ~~ c^1=\fr {y^1}z, ~~  c^2=\fr {y^2}z,
  ~~
 t=\fr{y^1}{y^2} \equiv  \fr{c^1}{c^2}  \equiv  \fr ij,
\ee
assuming for definiteness that
$z>0, j>0,$ and $c^2>0.$


\s\s

\nin
{\bf Definition 3.9} The function $U=U(x,y)$ which determines the \FMF
 $F$ in accordance with the equality
$F= i_{\{3\}} U$ is called the {\it primary generating  function}.

\s  \s

\nin
Henceforth, we assume for the \fmfs s
$F=F(x,y)$ of the considered three-axes type
 the  particular  structure specified by the condition of
{\it separation of  angle dependence}
\be
  U=\bvU(x,\theta), ~ f=\bvf(x,\theta), ~ Z=\bvZ(x,\phi), ~ t=\bvt(x,\phi)
\ee
of the {\it characteristic functions}
 $\{U,f,Z,t\}$
which enter the  functions $F$ in accordance with the following
{\it \SScht}:

\s

\be
F=z  \hU(x,f), ~~ f=c^2  \hZ(x,t).
\ee


It follows that
\be
  \bvU(x,\theta)=\hU(x,\bvf(x,\theta)), ~~ \bvZ(x,\phi)=\hZ(x,\bvt(x,\phi)),
\ee
and

\be
 \theta=\breve\theta(x,f), ~~
\phi=\breve\phi(x,t).
\ee


Since the functions $U$ and $f$ are homogeneous of degree zero with respect to  tangent vectors $\{y^i\}$,
it is also convenient to use the \reprns
\be
U=\ccUss(x,c^1,c^2), ~~ f=\ccfss(x,c^1,c^2),
\ee
where
$  \ccUss=\hU(x,\ccfss).
$


\s  \s

\nin
{\bf Definition 3.10} A Finsler space is called the
{\it \UZfsp}
if the \fmf $F$ of the space obeys the
\SSch
 (3.9). Such metric functions
will be denoted symbolically by
\FUZs.

\s\s


For the components of the unit vector $l^i=y^i/F$ we obtain the following angle representations:

\be
l^3=\fr1{\bvU(x,\theta)}, ~~  l^2=\fr{ \bvf(x,\theta)}{\bvU(x,\theta)}  \fr1{\bvZ(x,\phi)},
~~
l^1=l^2\bvt(x,\phi)\equiv
\fr{ \bvf(x,\theta)}{\bvU(x,\theta)}  \fr{\bvt(x,\phi)}{\bvZ(x,\phi)}.
\ee
It is instructive to compare the introduced functions  with their Euclidean precursors:
\s

$$
 \bvU^{Euclidean}=\fr1{\cos\theta},  ~~  \bvf^{Euclidean}=\tan\theta, ~~ \bvZ^{Euclidean}=\fr1{\sin\phi},
 ~~ \bvt^{Euclidean}=\fr1{\tan\phi},
$$
and
$$
l^{3;Euclidean}=\cos\theta, ~~  l^{2;Euclidean}=\sin\theta\sin\phi,
~~
l^{1;Euclidean}=\sin\theta \cos\phi.
$$
It can be said that our method consists in extending the Euclidean trigonometric functions
to obtain the generalized {\it \FRD trigonometric functions}:
\be
  \Cos(x,\theta)=\fr1{\bvU(x,\theta)}, ~~    \Sin(x,\theta)= \fr{\bvf(x,\theta)}{\bvU(x,\theta)},
   ~~ \Tan(x,\theta)=\bvf(x,\theta),
\ee
and

\be
   \Cos(x,\phi)=\fr{\bvt(x,\phi)}{\bvZ(x,\phi)}, ~~
   \Sin(x,\phi)=\fr{1} {\bvZ(x,\phi)}, ~~
   \Tan(x,\phi)=\fr1{\bvt(x,\phi)}, ~~
\ee
retaining the Euclidean structure of dependence on the angles $\theta$ and $\phi$:
\be
l^{3;\FRD}=\Cos(x,\theta), ~~  l^{2;Finsleroid}=\Sin(x,\theta)\Sin(x,\phi),
\ee
and
\be
l^{1;\FRD}=\Sin(x,\theta) \Cos(x,\phi).
\ee


With this trigonometric extension,  it is convenient to introduce
the {\it radius \RHx
of the  \hsec \FH} of the \frds:
\be
\RHxss  =\fr{\bvf}{\bvU}\Bigl|_{P_3}\equiv \Sin(x,\thtHxss),
\ee
where \thtHx corresponds to $P_3$.
This \RHx can conveniently be juxtaposed with the value
\be
\zHxss=\fr{1}{\bvU}\Bigl|_{P_3}  \equiv \Cos(x,\thtHxss).
\ee


The  conditions  (3.8)-(3.11) specify the dependence of the characteristic functions on
 tangent vectors $\{y^i\}$, and therefore the form of derivatives of the functions with respect to
$y^i$. In particular, we obtain the expansions
\be
l_i=z_iU+zU_i, ~~ U_i=\bvU_{\theta}\theta_i, ~~ \theta_i=\bvth_ff_i, ~~
\theta_{ij}=\bvth_{ff}f_if_j +  \bvth_ff_{ij},
\ee
{\it etc}, where
the subscripts $\theta$ and $f$ mean differentiations; $l_i=\partial F/\partial y^i, z_i=\partial z/\partial y^i,
U_i=\partial U/\partial y^i, f_i=\partial f/\partial y^i, f_{ij}=\partial f_i/\partial y^j, \theta_i=\partial \theta/\partial y^i,$ and
$\theta_{ij}=\partial \theta_i/\partial y^j$.

\s

Let us  use  the notation

$$
U_f=\D{ \hU}f, ~~  U_{ff}=\Dd{ \hU}{f}{f},~~ U_1=\D\ccUss{c^1}, ~~ U_2=\D\ccUss{c^2}, ~~
f_1=\D\ccfss{c^1}, ~~  f_2=\D\ccfss{c^2},
$$
such that

$$
\D f{y^3}=-\fr1zf, ~~ \D f{y^1}=\fr1zf_1, ~~ \D f{y^2}=\fr1zf_2.
$$


In many instances it is convenient to use the equalities
$U_1=U_ff_1$ and $U_2=U_ff_2$.
The identity
\be
f=c^1f_1+ c^2f_2
\ee
holds fine due to the  \hty
 implied.
Evaluating the covariant unit vector
components $ l_i=\partial  F/\partial{y^i}=l_1 i_i+l_2 j_i +l_3i_{\{3\}i}$
yields
\be
l_1=U_1, ~~ l_2=U_2, ~~ l_3=U-U_ff
\ee
[notice that $l_3=U-U_f( f_1c^1+ f_2c^2)$].


\s\s


With this preparation,  all the derivatives
$
l_{ij}=\partial l_{i}/\partial {y^j}
$
can  be found: at first, we arrive at
$$
y^3l_{33}=U_{ff}f^2, ~ y^3l_{31}=-U_{ff}ff_1, ~  y^3l_{32}=-U_{ff}ff_2,
$$
and after that we can use
$
U_1=U_ff_1
\equiv l_1
$
and conclude that
$
y^3l_{11}=
U_{ff}f_1f_1+U_ff_{11}, {\it etc.}
$
By following this method, we  find all the components
$
\hUZss_{ij}=Fl_{ij}
$
 of  the  tensor \hUZ and establish the validity of the following assertion.


\s\s

\nin
{\bf \Prpn 3.1} {\it The tensor
\hUZ
of the    \UZfsp
can be given by the following components:
\be
\hUZss_{33}=\cU_1f^2, ~~ \hUZss_{31}=-\cU_1ff_1, ~ ~ \hUZss_{32}=-\cU_1ff_2,
\ee
and
\be
\hUZss_{11}=
\cU_1f_1f_1+\cU_2f_{11}, ~~
\hUZss_{12}=
\cU_1f_1f_2+\cU_2f_{12},  ~~
\hUZss_{22}=
\cU_1f_2f_2+\cU_2f_{22},
\ee
where $\cU_1=UU_{ff}$ and $\cU_2=UU_{f}$.
}

\s\s

\nin
In turn, the explicit list of these representations determines explicitly all the components
  $\gUZss_{ij}=\hUZss_{ij}+l_il_j$  of  the \FMT \gUZ of the
  \UZfsps.

\s\s


Using the list, the attentive \calcns of the determinant
\be
\DUZss=\det(\gUZss_{ij})
 \ee
lead to the following significant result.

\s\s

\nin
{\bf Proposition 3.2} {\it  The determinant  $\DUZss$ of the \FMT  \gUZ of the
\UZfsp
is  the product
\be
\DUZss=D_1D_2
\ee
with the factors
\be
D_1=U^4U_{ff}U_f\fr1f, ~~  D_2=    Z^3   Z_{tt}
  \ee
which possess the property of separation of dependence on angles $\tht$ and $\phi$, namely
\be
D_1=D_1(x,\theta), ~~     D_2=D_2(x,\phi).
\ee
}


\s\s



\setcounter{equation}{0}

\nin
{\bf  \large 4  The \UZfsp  of \FRD type}

\s\s

\nin
The  representations (3.23) indicate  that subjecting $\theta$ to the nonlinear differential
equation
\be
\bvth_f\bvth_f=H^2\fr1UU_{ff}
\ee
is necessary and sufficient in order that  the \frd equalities
$$
h_{3c}=h^{\htph}_{3c}\equiv F^2\fr1{H^2}\theta_3\theta_c, ~~ c=1,2,3,
$$
hold, where
$\tht_c=\bvth_f\partial f/\partial y^c\equiv \partial\tht/\partial y^c$.
The vanishing $\prtl \phi/\partial  y^3=0$ has been taken into account.

\s\s

Thus we come to the equality

$$
h_{ij}=\fr1{H^2}\lf(\theta_i\theta_j+H^2c^2\fr1UU_fZ_{tt}t_it_j\rg)F^2,
$$
where $t_i=\partial t/\partial y^i$. Since $f=c^2Z$, we can write
\be
h_{ij}=\fr1{H^2}\lf(\theta_i\theta_j+H^2A_1A_2t_it_j\rg)F^2, ~~  A_1=\fr1UU_ff, ~ A_2= \fr1ZZ_{tt}.
\ee


\nin
This representation gives rise to
the validity of the following assertions.

\s  \s

\nin
{\bf Conclusive Theorem} {\it  The   \UZfsp
 is the \FRD space if and only if
the following two nonlinear differential equations are fulfilled:
\be
\bvth_f\bvth_f=H^2\fr1UU_{ff}, ~~ \fr1UU_ff=T\sin^2\theta,
\ee
with
\be
T=C\fr1{H^2},  ~~    T=T(x), ~~  C=C(x), ~~ H=H(x).
\ee
 At the same time, the condition of constancy of the indicatrix curvature
of the
\UZfsp
does not impose any equation on the dependence of the function $Z=Z(x,t)$ on the argument
$t$, as well as of the function $\bvZ=\bvZ(x,\phi)$ on $\phi$.
 To comply with the  constancy of the indicatrix curvature,
the  function  $\phi=\phi(x,t)$  must obey the equality
\be
\phi_t\phi_t=  C\fr {Z_{tt}}Z,
\ee
that is, the function must be given by the integral
\be
\phi=  \int\sqrt{C\fr {Z_{tt}}Z}\,dt.
\ee
The involved function $Z$ can depend on the argument $t$ in an arbitrary way.

}


\s\s

\nin
{\it Remark 4.1}   The  \UZfsp
 is the \FRD space at any choice of the dependence of $Z$ on $t$, provided
that the generating function $U$ fulfills the equations
(4.3).
This arbitrariness of the function $Z=Z(x,t)$ is rather unexpected phenomenon, and even a suspicious statement
for the first glance.
However, the explanation of the phenomenon can be formulated in simple words. Namely,
the condition of constancy of the indicatrix curvature of the   \UZfsp
 gives rise to the similar equations
$\bvth_f\bvth_f=H^2U_{ff}/U$ and $\phi_t\phi_t=  CZ_{tt}/Z$.
For the function $U$ the respective representation (2.7)  of the tensor $h^{\htph}_{ij}$ gives the second equation
$U_ff/U=T\sin^2\theta$
because of the presence of the function $\sin^2\tht$ in  the right-hand part of the representation.
No  functions of $\phi$ enters  the right-hand part, - whence no second equation can be obtained for the dependence of the function
$\bvZ$ on $\phi$.

\s\s


\nin
 {\it Remark 4.2} Owing to the assumed  Separation Scheme, the structure  of
dependence of the functions $A_1$ and $A_2$ on tangent vectors
is principally different, namely we have
$A_1=A_1(x,f)$ and $A_2=A_2(x,t)$ (see (4.2)).
Therefore, to get  the \frd representation
\be
\hUZss_{ij}=h^{\htph}_{ij}
\ee
with
$
h^{\htph}_{ij}=(1/H^2)\lf(\theta_i\theta_j+\sin^2\theta \phi_i\phi_j\rg)F^2
$
(see (2.7)), where $\phi=\phi(x,t)$ and $\phi_i=\phi_tt_i$,
the  equations
$
(1/U)U_ff=T\sin^2\theta$ and $ \phi_t\phi_t= C (1/Z)Z_{tt}
$ with    $T=T(x), ~  C=C(x)$, and $ T=C/H^2$
must be fulfilled.

\s\s

\nin
{\bf Definition 4.1}
The equations (4.3)-(4.4)  formulated in the above Conclusive Theorem
will be called the
{\it \UZ-Finsleroid characteristic equations},
the \it {\FE-equations} \rm
for short.

\s\s


\nin
{\bf Proposition 4.1}
{\it
The dependence of the characteristic functions $U=\bvU(x,\tht)$ and  $f=\bvf(x,\tht)$
of the \UZ-Finsleroid space
on the argument $\tht$
can explicitly be found upon integrating  the \FE-equations  \rm(4.3)-(4.4), \it namely the respective representations are given by  Formulas \rm(4.19) \it and
 \rm (4.22).

\s\s


Let us arrive at the solutions.
To this end we differentiate the second \FE-equation
$(1/U)U_ff=T\sin^2\theta$
with respect to $f$, whereupon using the first
\FE-equation $\bvth_f\bvth_f=H^2(1/U)U_{ff}$,
which yields
$$
-\fr1f  T\sin^2\theta T\sin^2\theta +
\fr1{H^2}(\bvth_f)^2 f+
 \fr1f T\sin^2\theta
 =
  2T\bvth_f\sin\theta \cos\theta.
   $$
This equation can be reduced to

$$
-T\sin^4\theta +
\fr1{TH^2}  (\bvth_ff)^2 +
\sin^2\theta(\sin^2\theta+\cos^2\theta)
 =
  2f\bvth_f\sin\theta \cos\theta.
   $$
Noting  $TH^2=C$, we arrive at the quadratic equation
$$
\lf(\fr1{\sqrt C}\bvth_ff-\sqrt C\sin\theta \cos\theta\rg)^2
 =
(T -1)\sin^4\theta
+(C-1)  \sin^2\theta\cos^2\theta.
   $$
Inserting
\be
\bvth_ff =
 C
 R_{17} \sin\theta,
\ee
we obtain the equation
$$
C\lf(R_{17}- \cos\theta\rg)^2
 =
(T -1)\sin^2\theta
+(C-1)  \cos^2\theta
   $$
from which the function $R_{17}$ can be found, namely
\be
R_{17}=
 \cos\theta
+
\sqrt{L_{17}},
\ee
where
\be
L_{17}=
\lf(\fr1{H^2} -\fr1C\rg)\sin^2\theta
+\lf(1 -\fr1C\rg)  \cos^2\theta.
\ee

\s


\s

The last function can conveniently be written in the form
\be
 L_{17}=P_1+H_1  \sin^2\theta \equiv H_2-H_1\cos^2\theta,
\ee
 where
\be
P_1=1-\fr1{C}, ~~ H_1=\fr1{H^2}-1, ~~ H_2=\fr1{H^2}-\fr1{C}.
\ee
The argument dependence of this function is of the type
$ L_{17}= L_{17}(x,\tht).$

\s

In order to find the dependence of the functions $U$ and $f$ on the argument $\theta$
we introduce the functions
\be
 I=
\exp\lf(-\sqrt{H_1}
\arcsin \lf(\sqrt{\fr{H_1}   {H_2}}
\cos\theta\rg)\rg)
\ee
and

\be
 Y=  \lf(\fr1{\sqrt{H_2}\sin\theta}(\sqrt{ L_{17}}+\sqrt {P_1}\cos\theta) \rg)^2
 \equiv
\fr{\sqrt{P_1}\cos\theta+\sqrt{ L_{17}}} {-\sqrt{P_1} \cos\theta+\sqrt{ L_{17}}},
\ee
together with
\be
 Y= \fr{\cY_+}{\cY_-} \equiv \lf(  \fr{\cY_+}{\sin\theta}\rg)^2,
\ee
where
\be
\cY_+=\sqrt{\sin^2\theta+a\cos^2\theta}+\sqrt a\cos\theta, ~~
\cY_-=\sqrt{\sin^2\theta+a\cos^2\theta}-\sqrt a\cos\theta
\ee
with $ a=P_1/H_2.$
The function
\be
 Y_1= Y^{\frac12\sqrt{P_1}}
\ee
enters various significant \frd representations.

\s



With

$$
\fr1UU_ff=T\sin^2\theta,  ~~  \bvth_ff =  CR_{17} \sin\theta , ~~  T=C\fr1{H^2},
$$
we come to the following differential equation for the generating  function $U=\bvU$:

\be
 \fr1{\bvU}\bvU_{\theta}
=
\fr1{H^2R_{17}}\sin\theta.
\ee
Integrating yields the   explicit representation
\be
\bvU=
 C_{22}\fr1{R_{17}}I, ~~    C_{22}=C_{22}(x),
\ee
which can readily be verified.

\s


The last formulas entail
\be
 \lf(\fr1{\bvU}\rg)_{\theta}=-\fr1{ C_{22}} \fr1{H^2}\fr1I \sin\theta.
\ee


 Also,
we should  solve the equation (4.8) to determine  the dependence of  $f=\bvf(x,\theta)$ on $\theta$.
To this end we write the equation as
\be
(\ln \bvf)_{\theta}=
\fr1{
CR_{17} \sin\theta
}.
\ee
On the other hand,
$$
\lf(\ln \lf(Y_1\fr1{R_{17}}\sin\theta \rg)\rg)_{\theta}=
\fr1{CR_{17}\sin\theta}.
$$
We observe that the equality
$$
(\ln \bvf)_{\theta}=
\lf(\ln \lf( \fr{\sin\theta}{R_{17}}Y_1\rg)\rg)_{\theta}
$$
holds. Therefore, the sought solution reads
\be
\bvf=
 C_{33}   \fr{\sin\theta}{R_{17}}Y_1, ~~  C_{33}=C_{33}(x).
 \ee
Thus   \Prpn 4.1 is valid.

\s\s

With the help of  (4.21) we can obtain
\be
 \bvf_{\theta}=
 C_{33}\fr1{
C(R_{17})^2}
Y_1.
\ee

It follows also that

\be
\fr1{\bvU}
\bvf=
\fr{ C_{33} }{ C_{22}}  \fr1I Y_1\sin\theta
 \ee
is the radius of the \hsec (cf. (3.18)).


The \FRD trigonometric functions (see (3.14)) admit the following explicit
representations:

\be
\Sin(x,\tht)=
\fr{ C_{33} }{ C_{22}}  \fr1I Y_1\sin\theta,  ~~
\Cos(x,\tht)=
\fr1{ C_{22}} R_{17}\fr1I.
  \ee


If we compare two formulas
$ (\bvth_ff)^2=H^2(1/U)U_{ff}f^2$ and $ \bvth_ff = C  R_{17} \sin\theta $
(see (4.3)  and (4.8)),
we can obtain the equality
\be
\fr1UU_{ff}f^2=\fr1{H^2}    C^2 ( R_{17})^2 \sin^2\theta.
\ee
Owing to
$ U=  C_{22}(1/R_{17})I$
(see  (4.19)),            
we can also propose the representation
\be
UU_{ff}=\fr1{H^2}    C^2 (C_{22})^2I^2 \sin^2\theta
\ee
with very simple dependence on the variable $\tht$ in the
\rhpt.

Taking into account the equalities
$$
  f=  C_{33}   \fr{\sin\theta}{R_{17}}Y_1, ~~   U=  C_{22}\fr1{R_{17}}I
$$
(see (4.19) and (4.22)), we can conclude that
$$
\fr1UU_{ff}=\fr1{H^2(C_{33}  )^2(Y_1)^2}    C^2 ( R_{17})^4
$$
and
\be
U^3U_{ff}=\fr{C^2 (C_{22})^4}{H^2(C_{33}  )^2}    \fr1{(Y_1)^2}I^4.
\ee


\s\s

The representations (4.26)-(4.28)
 help us to evaluate the determinant
$ D=\det(g_{ij})$
of the   \frdmt
 $g_{ij}=g_{ij}(x,y)$.
Namely,
 we can elucidate the dependence of the factor $D_1$ on the argument $\tht$ in the determinant
 $D=\det(g_{ij})=D_1D_2 $ (see (3.26)).
To this end, the representations
$$
\fr1UU_{ff}f^2=\fr1{H^2}    C^2 ( R_{17})^2 \sin^2\theta , ~~   \fr1UU_ff=T\sin^2\theta,
$$
(see (4.26) and (4.3))
can conveniently be applied in
 (3.27), which yields
$$
D_1=U^6\fr1{H^2}    C^2 ( R_{17})^2 \sin^2\theta
T\sin^2\theta\fr1{f^4}.
$$
Inserting here the functions
$$
 U=  C_{22}\fr1{R_{17}}I, ~~  f=
 C_{33}   \fr{\sin\theta}{R_{17}}Y_1
$$
(see (4.19) and (4.22))
we get the following result:
\be
D_1=C_{11}
 I^6\fr1{(Y_1)^4}, ~~    C_{11}=C_{11}(x), ~~
 C_{11}=\fr1{H^2}    C^2T \fr{(C_{22})^6}{( C_{33})^4}.
\ee





\s

\nin
{\bf  \large 5 Conclusions}

\s\s

\nin
The established \CT yields the exhaustive answer to the question ``What is the form of the ODE-s which describe
 the azimuthal  $\tht$-angle dependence of the   three-axes \frd {\chfns}s provided that the
\SSch
is used?''

At the same time, the derived ODE-s
don't impose any restriction on the dependence of the {\chfns}s on the \frd  polar angle $\phi$. The dependence
can be introduced or specified by means of additional geometrical conditions.

The clarification of the regularity properties of the  obtained solutions  requires the systematic and attentive
study.

\s\s

\def\bibit[#1]#2\par{\rm\noindent\parskip1pt
                     \parbox[t]{.05\textwidth}{\mbox{}\hfill[#1]}\hfill
                     \parbox[t]{.925\textwidth}{\baselineskip11pt#2}\par}

\nin {\bf \large References}

\s\s

\bibit[1]  Asanov,G.S.:Pseudo-Finsleroid metric function of spatially anisotropic relativistic type
 (2015).  arXiv:1512.02268


\bibit[2]  Asanov,G.S.:Pseudo-Finsleroid metrics with two axes,  European Journal of Mathematics.  {\bf 3}(4), 1076-1097  (2017).
 DOI 10.1007/s40879-017-0160-6

\s

\bibit[3] Asanov,G.S.:Two-axes pseudo-Finsleroid metrics: general overview and angle-regular solution  (2017).
 arXiv:1709.02683


\bibit[4] Cartan,E.:Les Espaces de Finsler.   Actualit\'es Scientifiques et Industrielles, vol. 79.
 Hermann,  Paris  (1934)

\s

\bibit[5]  Berwald,L.:\"Uber Finslerische und verwandte R\"aume,  Cas. Mat. Fys.  {\bf 64}, 1-16 (1935)


\bibit[6]  Busemann,L.:The Geometry of Finsler Spaces,  Bull. Amer. Math. Soc.  {\bf 56}, 5-15 (1950)


\bibit[7]   Rund,H.:The Differential Geometry of Finsler  Spaces. Die Grundlehren der Mathematischen Wissenschaften, Band, vol. 101.  Springer, Berlin, 1959


\bibit[8] Horvath,J.I.:New geometrical methods of the theory of physical fields, Nuov. Cim. 9 (1958)  444-496.


\bibit[9]    Ingarden,R.S.:On physical applications of Finsler Geometry, Contemporary Mathematics 196 (1996) 213-223.


\bibit[10]  Asanov,G.S.:Finsler Geometry, Relativity and Gauge  Theories,  D.~Reidel Publ. Comp., Dordrecht 1985.









\bibit[11] Bao,D.,Chern,S.-S.,Shen,Z.:An Introduction to Riemann-Finsler Geometry,
Graduate Texts in Mathematics, vol. 200.  Springer, New York  (2000)









\bibit[12]  Matveev,V.S.,Papadopoulos,A.,Rademacher,H.-B.,Sabau,S.V. (Eds.):
The Topical Issue "Finsler Geometry: New Methods and Perspectives"
  European Journal of Mathematics (20017)3: 763-766.
DOI 10.1007/s40879-017-0195-8


\bibit[13]  European Journal of Mathematics,  December 2017, Issue 4, Pages 763-1273.
Special Issue: Finsler Geometry, New Methods and Perspectives





\bibit[14]   Asanov,G.S.:Finsler connection properties generated by the two-vector angle developed on the indicatrix-inhomogeneous  level, Publ. Math. Debrecen {\bf 82}(1) (2013) 125-155


\bibit[15] Vincze,Cs.:On Asanov's Finsleroid-Finsler metrics as the solutions of a
conformal rigidity problem,
 Differential Geometry and its Applications  {\bf 53} 148-168    (2017)

\end{document}